\documentclass[a4paper,12pt]{article}

\usepackage{amsmath,amssymb,amscd,theorem}
\usepackage{xypic}

\usepackage[english]{babel}

\newtheorem{theorem}{Theorem}[section]
\newtheorem{proposition}[theorem]{Proposition}
\newtheorem{lemma}[theorem]{Lemma}
\newtheorem{corollary}[theorem]{Corollary}

{\theorembodyfont{\rmfamily }
\newtheorem{definition}[theorem]{Definition}

\newtheorem{remark}[theorem]{Remark}
\newtheorem{example}[theorem]{Example}
\newtheorem{examples}[theorem]{Example}}

\def\RR{\mathbb{R}}
\def\CC{\mathbb{C}}
\def\d{\mathrm{d}}
\def\proof{\noindent \textit{Proof: }}
\def\qed{\hfill $\square$ }
\def\Isom{\mathrm{Bihol}(X)}

\begin{document}

\title{NATURAL  OPERATIONS \\ ON HOLOMORPHIC FORMS}

\date{}

\author{Navarro, A.\thanks{Institut f\"{u}r Mathematik, Universit\"{a}t Z\"{u}rich, Switzerland.} 
\and Navarro, J.\thanks{Corresponding author. Email address: \texttt{navarrogarmendia@unex.es},
\newline Department of Mathematics, Universidad de Extremadura,  Spain.}
\and Tejero, C.\thanks{Department of Mathematics and IUFFyM, Universidad de Salamanca, Spain.}
}

\maketitle

\begin{abstract}

We prove that the only natural differential operations between holomorphic forms on a complex manifold are those obtained using linear combinations, the exterior product and the exterior differential. In order to accomplish this task we first develop the basics of the theory of natural holomorphic bundles over a fixed manifold, making explicit its Galoisian structure by proving a categorical equivalence {\it \`{a} la Galois}.

\medskip


\end{abstract}

\tableofcontents

\section*{Introduction}

Roughly speaking, the term {\it natural operation} makes reference to those operations in differential 
geometry whose definition does not depend on choices of coordinates. As an example, if $M$ is a smooth manifold, the exterior differential
$$ \d \colon \Lambda^p T^*M \, \rightsquigarrow \, \Lambda^{p+1} T^*M $$ is an $\RR$-linear differential operator whose definition
can be made on a chart, but it is independent on the chosen coordinates.

The study of these constructions is already present in the very beginnings of
differential geometry, although its modern and systematic development was established
in the 1970s, due to the works of Nijenhuis (\cite{Nijenhuis_72}), Atiyah-Bott-Patodi
(\cite{Atiyah_73}), Epstein-Thurston (\cite{Thurston_76}) and Terng
(\cite{Terng_78}), among many others. In a sense, this development culminated in the
monograph by Kol\'{a}\v{r}-Michor-Slov\'{a}k (\cite{KMS}), which has become the standard
reference in the subject since then.

In this paper, we are interested in the following characterization of the exterior
differential, due to Palais:

\medskip
\noindent {\bf Theorem (\cite{Palais_59}):} {\it Let $\,M\,$ be a smooth manifold. If
$$\,P \colon \Lambda^p T^*M
\rightsquigarrow \Lambda^{p+1} T^*M \,$$ is a natural and $\RR$-linear differential
operator, then there exists $\,\lambda \in \RR\,$ such that $\,P \, =
\, \lambda \, \d$.}
\medskip

This result was later improved by Krupka-Mikol\'{a}\v{s}ov\'{a} (\cite{Krupka_84}), who
dropped the hypothesis of linearity, and considered (non-linear) differential operators of
arbitrary order $\leq r$:

\medskip
\noindent {\bf Theorem (\cite{Krupka_84}):} {\it Let $\,M\,$ be a smooth manifold and let $p> 0$. If
$$\,P \colon \Lambda^p T^*M
\rightsquigarrow \Lambda^{p+1} T^*M \,$$ is a natural differential operator of order
$\leq r$, then there exists $\,\lambda \in \RR\,$ such that $\,P \, = \, \lambda\, \d$.}
\medskip

The case $p=0$ lies outside this characterization; as an example, the differential
operator $\, \mathcal{C}^\infty (M) \rightsquigarrow T^*M \, , \, f \mapsto f \, \d f
\,$ is natural and it is not a constant multiple of the exterior differential.

In their monograph, Kol\'{a}\v{r}-Michor-Slov\'{a}k (\cite{KMS}, Prop. 25.4) refined this statement 
even more, by considering regular (and natural) local operators. Theirs is a two-step proof: firstly, 
they proved that regular local operators are indeed differential operators--the Peetre-Slov\'{a}k's
theorem--and then, a version of the result above applies.

Quite recently, there have been a couple of reformulations of this result
(\cite{Katsylo_08}, \cite{Navarro_2015}) that, essentially, confirm that
the only natural operation between differential forms is the exterior differential.
These theorems have been used for several different applications (\cite{FreedHopkins}, 
\cite{Katsylo_08}, \cite{Mason_2014}, \cite{Navarro_2015}, \cite{Timashev_15}). As an example, these ideas have been succesfully extended to study differential operations on forms on contact manifolds, allowing a deeper understanding of the Rumin operator (\cite{Bernig}).

In this paper, we aim to study up to what extent this statement remains true for holomorphic forms over a complex analytic manifold. 

Katsylo and Timashev (\cite{Katsylo_08}, Theorem 3.1) have already considered natural, $\mathbb{C}$-linear differential operators, obtaining--among other things--a result analogous to the aforementioned theorem of Palais. 

Our main result goes one step further and removes the hypothesis of linearity (see also Theorem \ref{General}):

\medskip
\noindent {\bf Theorem \ref{Main}:} {\it Let $X$ be a complex analytic manifold of
dimension $n$.

Let $\,\mathbb{C}\{ u,v\}\,$ be the anti-commutative algebra of polynomials on the
variables $\,u,v$, of degree $\,\mathrm{deg}\, u=p$, $\,\mathrm{deg}\, v=p+1$ (Definition \ref{AntiConmutativos}).

Any natural differential operator of order $\leq r$
$$ P\colon \Lambda^{p} T^*X \  \rightsquigarrow\ \Lambda^q T^*X \ ,  $$
where $0\leq q\leq n= \dim X$, can be written as
$$P\left(j^r_x \omega \right)\,=\, \mathbf{P}(\omega,\mathrm{d}\omega )_x \quad , $$
for a unique homogeneous polynomial $\,\mathbf{P}(u,v)\in\mathbb{C}\{ u ,v \}\,$ of
degree $\,q$.}
\medskip

The proof closely follows the one given in \cite{Navarro_2015} for the case of differentiable manifolds, which, in turn, applies
techniques already used by the Czech school (\cite{KMS}, \cite{Krupka_84}). 

Nevertheless, in order to accomplish our task it is necessary to deal with natural bundles in the category of complex manifolds and holomorphic maps. Among experts, it is well-known that a great part of the theory in the smooth case  is purely algebraic, so that it extends without difficulties to this setting. This is the approach followed by Katsylo and Timashev (\cite{Katsylo_08}), although their definition of $k$-th order natural bundles over a manifold as associated bundles to the $k$-th order coframe bundle is not conceptually the most fundamental one. We will see that it assumes, as a starting point, a key result of the theory: the theorem of equivalence (\ref{TeoremaEquivalencia}). On the other hand, the usual approach to natural bundles (\cite{KMS}) as functors on the category of manifolds and local diffeomorphisms is the most general possible. However in this way one is not able to capture an important internal property, the Galoisian structure, enjoyed by natural bundles over a fixed manifold.

Thus,  we devote the first section of this paper to developing the basics of  the theory on natural holomorphic bundles over a fixed manifold from the most fundamental definitions and at the same time to make transparent its natural Galoisian structure. This will also serve to fix the notation used in the rest of the paper. We follow an approach that clearly resembles the modern expositions of Galois theory, where the role of the Galois group is played here by the group of finite order jets of biholomorphisms leaving a point fixed.  
Our exposition is for the most part rather schematic and we just prove several key points, since most of the ideas carry along from 
the usual presentations in the differentiable case. 

However, we have identified several important results in the smooth theory 
whose holomorphic counterparts do not follow straightforwardly and therefore we do not know yet wether they do hold or not in the holomorphic category. In our opinion these problems deserve further attention in order to have a complete understanding of the theory of holomorphic natural bundles and its relation with its smooth counterpart. We plan to address these questions in the future but now let us
just point out three of them:

\begin{itemize}
\item Is every holomorphic natural bundle a finite order one?

\item Can the regularity condition in the definition of natural bundle be deduced from the other two?

\item Does Peetre-Slov\'{a}k's theorem on the finiteness of the order of regular local operators hold in the holomorphic category? 
\end{itemize}

Probably, the solution of these open problems would pave the road to put this holomorphic theory on an equal footing 
to the smooth, classical one.

\section{Holomorphic natural bundles}

Let us present in this section a non-standard definition of natural bundle, first proposed in the smooth case by Juan B. Sancho Guimer\'{a}. We use it due to its simplicity and its benefits in order to state a Galois-type theorem. Later, in Section \ref{AlternativeDefi}, we will explain the equivalence with a functorial exposition that is analogous to the one commonly used in the smooth theory.

\medskip

Let us fix a complex analytic manifold $X$ of dimension $n$.

If $\,\pi\colon F \to X\,$ is a holomorphic fibre bundle (i.e., a holomorphic map isomorphic, locally, on $X$,  to the trivial fibration $U\times S \to U$), then a {\it lift to $F$} of an analytic isomorphism $\,\tau \colon U \to V\,$ between open sets in $X$ is any
analytic isomorphism $\,\ell^F(\tau) \colon F_{U} := \pi^{-1}(U) \to F_{V} := \pi^{-1}(U)\,$ such that the following square is commutative:
$$
\xymatrix{
F_{U} \ar[r]^-{\ell^F(\tau)} \ar[d]_-{\pi}  & F_{V} \ar[d]^-{\pi} \\
U \ar[r]^-{\tau} &  V} 
$$

If $Z$ is a complex manifold we denote by  $\,\mathrm{Bihol}(Z) \, $ the set of analytic isomorphisms
$\tau \colon U \to V$ between open sets in $\,Z$.

\begin{definition}\label{DefinitionJS}
A natural bundle over $X$ is a pair formed by a holomorphic fibre bundle $F \to X$, together with a lifting of biholomorphisms:
$$\ell^F\colon \Isom  \to \textrm{Bihol}(F) \quad , \quad \tau \ \mapsto \ \ell^F(\tau) $$ that
satisfies:
\begin{itemize}
\item {\it Functoriality:} $\ell^F(\mathrm{Id}) = \mathrm{Id}$ and $\ell^F(\tau \circ \tau')=\ell^F( \tau)\circ \ell^F(\tau')$.
\item {\it Locality:} For any analytic isomorphism $\tau \colon U \to V$ between open sets on $X$, and for any open set $U' \subset U$:
$$
\ell^F(\tau_{|U'}) = \ell^F(\tau)_{|F_{U'}} \ .
$$
\item {\it Regularity:}  If $\{ \tau_t \colon U_t \to V_t \}_{t\in T}$ is an analytic family of biholomorphisms between open sets in $X$, parametrized by an analytic manifold $\,T$, then the collection of lifts $\,\{ \ell^F(\tau_{t} )\colon F_{U_t} \to F_{V_t} \}_{t \in T}\,$ is an analytic family of biholomorphisms between open sets in $\,F$.
\end{itemize}

Let $F, \bar{F} \to X$ be natural bundles over $X$. A morphism of natural bundles (or natural morphism) on an open set $W$ 
is a morphism of bundles  $\varphi \colon F_W \to \bar{F}_W\, $ that commutes with the action of $\,\Isom \,$; that is, such that
for any  analytic isomorphism $\tau \colon U \to V$ between open sets on $W$, the following square is commutative:
$$
\xymatrix{
F_{U} \ar[r]^-{\varphi} \ar[d]_-{\ell^F(\tau) }  & \bar{F}_{U} \ar[d]^-{\ell^{\bar F}(\tau)} \\
F_{V} \ar[r]^-{\varphi}  & \bar{F}_{V}
} 
$$

A natural bundle $F \to X$ has order $k$ if for any biholomorphisms $\tau,\tau'\colon U\to V$ defined between open sets of $X$, and any $x\in U$, it holds
$$j^k_x\tau=j^k_x\tau'\quad\Rightarrow\quad \ell^F(\tau)=\ell^F(\tau')\quad\text{when restricted to the fibre }F_x \ .$$

The category of natural bundles of order $ k$ over $X$ will be denoted $\textbf{Nat}_X^k$.
\end{definition}

\begin{remark} From now on, if ($F\to X,\ell^F)$ is a natural bundle and $\tau \colon U \to V$ is a biholomorphism between open sets of $X$, whenever there is no risk of confusion we will denote the lift $\ell^F(\tau)\colon F_U\to F_V$ simply by $\tau_*\colon F_U\to F_V$.
\end{remark}

\begin{examples}
The trivial bundle $F = Y \times X \to X$, with the lifting $\tau_* (y, x) = \left( y , \tau (x) \right) $, is a natural bundle of order 0.

The tangent bundle $TX\to X$ is a natural bundle: the lifting of an analytic isomorphism $\tau\colon U\to V$ is its tangent linear map $\tau_*\colon TU\to TV$. More generally, any tensor bundle $\otimes^pT^*X\otimes^qTX$ is natural, of order $1$.

If $F\to X$ is a natural bundle of order $k$, then its bundle of $r$-jets, $\,J^rF\to X$, is also a natural bundle, of order $k+r$.
\end{examples}

\begin{example}[Universal natural bundle of order $ k$]\label{FibradoUniversal} Let us fix a po\-int $\, x_0 \in X$.

It is straightforward to check that the set $\mathcal{U}^k_{x_0}(X)$ formed by the jets  $j^k_{x_0}\sigma$ of germs of analytic isomorphisms $\,\sigma\colon U\rightarrow V\,$ between open sets of $X$ such that $x_0\in U$, has a natural structure of complex manifold that makes the projection:
$$\begin{CD}
\pi \colon \mathcal{U}^k_{x_0} (X)@>>> X\qquad ,\qquad j^k_{x_0}\sigma\mapsto \sigma(x_0)\ .
\end{CD}$$
into a holomorphic locally trivial bundle that is called the universal bundle of order $k$. 

This bundle is natural of order $ k$. The lift of a biholomorphism  $\tau\colon U\to V$ between open sets of $X$ is defined as:
$$\begin{CD}
\mathcal{U}^k_{x_0}(X)_{|U} @>{\tau_*}>> \mathcal{U}^k_{x_0}(X)_{|V}\quad ,\qquad \tau_*(j^k_{x_0}\sigma):=j^k_{x_0}(\tau\circ\sigma)\ .
\end{CD}$$

Moreover, let $\,G^k_{x_0}\,$ be the complex Lie group of $k$-jets $\,j^k_{x_0}\xi\,$ of germs of biholomorphisms $\xi\colon U\rightarrow U'$  leaving the point $x_0$ fixed. This group acts on the right on $\mathcal{U}^k_{x_0}(X)$:
$$ (j^k_{x_0}\sigma) \cdot(j^k_{x_0}\xi):=j^k_{x_0}(\sigma\circ  \xi). $$ Via this action, the universal bundle becomes a principal $G^k_{x_0}$-bundle.
\end{example}

\subsection{Theorem of equivalence}

Let $X$ be a complex analytic manifold, and let $x_0 \in X$ be any point. Recall that $\,G^k_{x_0}\,$ denotes the complex Lie group of $k$-jets $\,j^k_{x_0}\xi\,$ of germs of analytic isomorphisms $\xi \colon U \rightarrow  U'\,$  leaving the point $x_0$ fixed.

\subsubsection*{Fibre functor}

The condition of a natural bundle $F \to X$ being of order $\leq k$ amounts to saying that $G^k_{x_0}$ acts on the left on the fibre $F_{x_0}$ over the point $x_0 \in X$: for any point $\,e \in F_{x_0}\,$ and any $\,g= j^k_{x_0} \xi \in G^k_{x_0}$, define
\begin{equation}\label{Action}
g \cdot e := \xi_* (e) \ .
\end{equation}

This action is holomorphic: in a chart around $\,x_0$, it is easy to check the existence of an analytic family of isomorphisms $\,\tau_g \colon U_g \to V_g$, parametrized by $\, g\in G^k_{x_0}$, such that $j^k_{x_0} \tau_g = g $. By the regularity hypothesis, $(\tau_g)_*$ is an analytic family of isomorphisms of $F_{x_0}$; hence the action defined in (\ref{Action}) is holomorphic in $g$.

Moreover, if $F$, $\bar{F}$ are natural bundles of order $ k$, then any morphism of natural bundles $\varphi \colon F \to \bar{F}$ obviously defines a $G^k_{x_0}$-equivariant analytic map $\varphi_{x_0} \colon F_{x_0} \to \bar{F}_{x_0}$.

Thus, if $G^k_{x_0}\textbf{-Man}$ denotes the category of analytic $G^k_{x_0}$-manifolds, we have constructed a ``fibre functor'':
$$
\mathcal V_{x_0}\colon \textbf{Nat}_X^k  \xrightarrow{\ \  \ \ } G^k_{x_0}\textbf{-Man}    \quad , \quad F \ \mapsto \ F_{x_0}  \ .
$$

\subsubsection*{Associated bundle functor}

Let $F_0$ be a $G^k_{x_0}$-manifold; that is, $F_0$ is a complex manifold together with an analytic left action:
$$\begin{CD}
 G^k_{x_0}\times F_0 @>{\cdot}>> F_0  \ .
 \end{CD}$$

Let $\,\mathcal{U}^k_{x_0}(X)\, $ be the universal natural bundle of order $k$ (Example \ref{FibradoUniversal}). The group $G^k_{x_0}$ acts on the right on $\,\mathcal{U}^k_{x_0}(X)\times F_0\,$: $$ (j^k_{x_0}\sigma,e) \cdot(j^k_{x_0}\xi):=((j^k_{x_0}\sigma) \cdot(j^k_{x_0} \xi), (j^k_{x_0}\xi)^{-1}\cdot(e)). $$

\begin{definition}
For any $G^k_{x_0}$-manifold  $F_0$, its associated bundle is defined as:
$$\begin{CD}
F := (\mathcal{U}^k_{x_0}(X) \times F_0)/G^k_{x_0} @>>> X\quad ,\qquad [(j^k_{x_0}\sigma, e)]\mapsto \sigma(x_0) \ .
\end{CD}$$
\end{definition}

The associated bundle $F\to X$ is natural of order $\leq k$: the lifting of a biholomorphism $\tau\colon U\to V$ between open sets of  $X$ is defined by the formula:
 $$\begin{CD}
 F_U=(\mathcal{U}(X)^k_{\, U}\times F_0)/G^k_{x_0} @>{\tau_*}>> (\mathcal{U}(X)^k_{\, V}\times F_0)/G^k_{x_0}=F_V
  \end{CD}$$
  $$\tau_*\left([j^k_x\sigma,e]\right):=[\tau_*(j^k_x\sigma),e] \ . $$

This lifting is well defined because the actions of $\tau_*$ and $G^k_{x_0}$ on $\mathcal{U}^k_{x_0}(X)$ commute:
\begin{align*}
\tau_*\left((j^k_{x_0}\sigma)\cdot(j^k_{x_0}\xi)\right) \, &= \, \tau_*\left(j^k_{x_0}(\sigma\circ\xi)\right)
\, = \, j^k_{x_0}\left(\tau\circ\sigma\circ\xi\right) \\
&= \,\left(j^k_{x_0}(\tau\circ\sigma)\right)\cdot (j^k_{x_0}\xi) \, = \, \tau_*(j^k_{x_0}\sigma)\cdot(j^k_{x_0}\xi)\ .
\end{align*}

This construction is functorial: it transforms a morphism of $G^k_{x_0}$-manifolds $f \colon F_0 \to \bar{F_0}$ into the following morphism of natural bundles:
$$F=(\mathcal{U}^k_{x_0}(X)\times F_0)/ G^k_{x_0} \xrightarrow{\ \ \text{Id} \times f \ \ } (\mathcal{U}^k_{x_0}(X)\times \bar{F_0}) / G^k_{x_0}=\bar{F} ,$$ that maps  $[j^k_x\sigma , e] $ to $\ [j^k_x\sigma , f(e) ] \ . $

Thus, we have constructed a functor:
$$
\mathcal A\colon G^k_{x_0}\textbf{-Man}  \xrightarrow{\ \  \ \ } \textbf{Nat}_X^k   \quad , \quad F_0 \ \mapsto \mathcal A(F_0) :=\ (\mathcal{U}^k_{x_0}(X) \times F_0) / G^k_{x_0}  \ .
$$

\subsubsection*{Theorem of equivalence}

\begin{theorem}\label{TeoremaEquivalencia}
Let $X$ be a complex analytic manifold and $x_0\in X$ be any point. Consider the complex Lie group $G^k_{x_0}$ of $k$-jets  of germs of analytic isomorphisms leaving the point $x_0$ fixed.

The ``associated bundle" functor $\mathcal A$ establishes an equivalence of categories:
$$\begin{CD}
 \mathcal A\colon  G^k_{x_0}\mathrm{\bf \text{-}Man} @= \mathrm{\bf Nat}_X^k  \ ,
\end{CD}$$
whose inverse functor is the fibre functor $\mathcal V_{x_0}\colon {\bf Nat}_X^k  \xrightarrow{\ \  \ \ } G^k_{x_0}{\bf \text{-}Man}$ .
\end{theorem}

\proof  Let us check that both functors are inverse to each other: if $F_0$ is a $G^k_{x_0}$-manifold and $F\to X$ is the associated bundle, then there is an isomorphism
$$\begin{CD}
F_0 @= \left[(\mathcal{U}^k_{x_0}(X)\times F_0)/G^k_{x_0}\right]_{x_0}=F_{x_0}\\
e & \mapsto & [(j^k_{x_0}\text{Id},e)]
\end{CD}$$
which is $G^k_{x_0}$-equivariant:
\begin{align*}
 j^k_{x_0}\xi\cdot [j^k_{x_0}\text{Id},e] &= \, \xi_*[j^k_{x_0}\text{Id},e]=
[\xi_*(j^k_{x_0}\text{Id}),e] \, = \, [j^k_{x_0}(\xi\circ\text{Id}),e] \\
& = \, [j^k_{x_0}\xi,e]=[j^k_{x_0}(\text{Id}\circ\xi),e)]=[j^k_{x_0}\text{Id},j^k_{x_0}\xi\cdot e] \ .
\end{align*}

Conversely, if $F \to X$ is a natural bundle of order $\leq k$, then the following map is a bundle isomorphism between $F$ and the bundle associated to $F_{x_0}$:
$$\begin{CD}
(\mathcal{U}^k_{x_0}(X)\times F_{x_0} ) / G^k_{x_0} @= F \\
[j^k_{x_0}\sigma , e_{x_0} ]  & \mapsto &  \sigma_* e_{x_0} \ .
\end{CD}$$

It is well defined: if we change the pair $(j^k_x\sigma , e_{x_0})$ by the equivalent element $(j^k_{x_0}\sigma , e_{x_0})\cdot \,j^k_{x_0}\xi=(j^k_{x_0}(\sigma\circ\xi),(j^k_{x_0}\xi)^{-1}\cdot e_{x_0})=(j^k_{x_0}(\sigma\circ\xi),\xi^{-1}_*e_{x_0})$, then:
$$[j^k_{x_0}(\sigma\circ\xi),\xi^{-1}_*e_{x_0}]\,\mapsto\, (\sigma\circ\xi)_*\xi^{-1}_* e_{x_0}=\sigma_* e_{x_0} \ .
$$

\qed

\begin{examples}
Natural vector bundles of order $ 1$ correspond to holomorphic linear representations of the complex general linear group $$G^1_{x_0}=\mathrm{GL}_{\mathbb{C}}(T_{x_0}X).$$ As an example, the natural bundle associated to the standard action of this group on $\otimes^pT_{x_0}^*X\otimes^qT_{x_0}X$ is the vector bundle of $(p,q)$-tensors.

Natural coverings of order $ k$ correspond to discrete $G^k_{x_0}$-manifolds. The complex group $G^k_{x_0}$ is connected, so that any (connected) holomorphic natural covering is trivial, in contrast with the smooth theory, where there exists the orientation covering (observe that any complex manifold is orientable).

If $\mathcal{O}_{x_0}$ denotes the ring of germs of holomorphic functions around $x_0 \in X$, then the group $G^k_{x_0}$ acts on $P^k_{x_0}:=\{ j^k_{x_0}f:f\in\mathcal{O}_{x_0}\}$ as follows: $\,(j^k_{x_0}g)\cdot(j^k_{x_0}f)=j^k_{x_0}(f\circ g^{-1})$. The corresponding natural bundle is the bundle $\, J^k(X,\mathbb{C})\to X\,$ of $k$-jets of holmorphic functions on $X$.
\end{examples}

\begin{corollary}\label{MorfismosNaturalesK} Let $F,\bar F\to X$ be natural bundles of order $k$.
Given any point $x_0 \,$ in an open set $U$, the assignment $\varphi\mapsto\varphi_{x_0}$ defines a bijection:
$$\begin{CD}
\left\{\begin{aligned}
&\mathrm{Morphisms\ of\ natural\ bundles \, } \\
& \hskip .7cm \ \, \ \varphi \colon F_U\to\bar F_U
\end{aligned}\right\} @=
\left\{\begin{aligned}
& G^k_{x_0}\text{--}\,\mathrm{equivariant\ analytic\ maps\, }\\
& \hskip 1cm \varphi_{x_0} \colon F_{x_0} \to\bar F_{x_0}
\end{aligned}\right\}
\end{CD} \ . $$

In particular, any morphism $\varphi$ of natural bundles is globally defined and covers the identity map on $X$ (i.e. $\varphi (F_y) \subseteq \bar{F}_y$, for any point $y \in X$). 
\end{corollary}

\begin{example} The only morphisms of  vector bundles $\varphi\colon TX\to TX$ that are natural are the homotheties, $\varphi(D)=\lambda D$, because:
$$\mathrm{Hom}_{\mathrm{nat}}(TX,TX)\, =\, \mathrm{Hom}_{G^1_{x_0}}(T_{x_0}X,T_{x_0}X)\, =\, \mathrm{Hom}_{\mathrm{GL}_{\mathbb{C}} (T_{x_0}X)}(T_{x_0}X,T_{x_0}X) \, = \,  \mathbb{C}\ .
$$
\end{example}

\subsection{Functorial definition}\label{AlternativeDefi}

Let $X, Y$ be complex manifolds of the same dimension $n$, and let $x_0 \in X$, $y_0 \in Y$ be arbitrary points. Any local isomorphism $\tau$ such that $\tau (x_0) = y_0$ induces an isomorphism $G^k_{x_0} \simeq G^k_{y_0}$. Hence, Theorem \ref{TeoremaEquivalencia} implies that a natural bundle over a complex manifold produces a natural bundle over any other complex manifold of the same dimension. This motivates a functorial approach to natural bundles, that we proceed to scketch in analogy with the standard references for smooth bundles (\cite{KMS}, \cite{Nijenhuis_72}). 

\medskip

Let ${\bf Man _{\CC}}[n]$ be the category of complex analityc manifolds of dimension $n$ and local analytic isomorphisms between them.

On the other hand, let ${\bf Bund}$ be the category whose objects are fibre bundles $F \to X$, and whose morphisms $$\,(f, \bar{f}) \colon \{ F\xrightarrow{\pi} X\} \longrightarrow \{\bar{F}\xrightarrow{\bar{\pi}} \bar{X}\}\,$$ are pairs of analytic maps $\,(f, \bar{f})\,$ making the following square commutative:
$$
\xymatrix{
F \ar[r]^-{f} \ar[d]_-{\pi }  & \bar{F} \ar[d]^-{\bar{\pi}} \\
X \ar[r]^-{\bar{f}}  & \bar{X}
} 
$$

\begin{definition}\label{StandardDefi}
A natural bundle in dimension $n$ is a covariant functor:
$$ \mathfrak{F} \colon {\bf Man_{\CC}}[n] \to {\bf Bund} $$ satisfying the following properties:
\begin{itemize}
\item If $\mathfrak{B} \colon {\bf Bund} \to {\bf Man_{\CC}}[n]$ denotes the base functor, $\mathfrak{B} (F \to X) := X$, then the composition $\mathfrak{B} \circ \mathfrak{F} \colon {\bf Man_{\CC}}[n] \to {\bf Man_{\CC}}[n]$ is the identity functor $\text{Id} _{{\bf Man_{\CC}}[n]} $.

For any complex $n$-manifold $X$, let us write $\,\mathfrak{F}(X)=\{\pi_X \colon \mathfrak{F}_X \to X \}$.

\item {\it Locality:} For any open inclusion $i\colon U \hookrightarrow  X$, the bundle $\mathfrak{F}_U$ is identified with $\pi_X^{-1} (U)$ via the map $\mathfrak{F}(i)$.

\item {\it Regularity:} If $\,\{ f_t\colon X_t\to Y_t\}_{t\in T}\,$ is an analytic family of local isomorphisms between $n$-manifolds, parametrized by a complex manifold $T$, then  $\,\{\mathfrak{F} (f_t) \colon\mathfrak{F}(X_t)\to\mathfrak{F}(Y_t)\}_{t\in T}\,$ is also an analytic family.
\end{itemize}

A morphism of natural bundles is a morphism of functors $\Phi \colon \mathfrak{F} \to \bar{\mathfrak{F}}$.

A natural bundle $\mathcal{F}$ has order $ k$ if for any local isomorphisms between $n$-manifolds $f,g \colon X \to Y$ and any point $x \in X$, the following holds:
$$ j^k_x f = j^k_x g \quad \Rightarrow \quad \mathfrak{F} (f)  = \mathfrak{F} (g) \  \mbox{ when restricted to the fibre } (\mathfrak{F}_X)_x \ .$$

Let us write $\mathrm{\bf Nat}^k_n$ fot the category of natural bundles in dimension $n$ of order $k$ just defined.
\end{definition}

\begin{proposition}\label{EquivalenciaDefiniciones} Let $\mathrm{\bf Nat}^k_n$ be the category of natural bundles in dimension $n$, of order $ k$ (Definition \ref{StandardDefi}).

Let $X$ be a complex analytic manifold of dimension $n$, and let $\mathrm{\bf Nat}^k_X$ denote the category of natural bundles over $X$, of order $\leq k$ (Definition \ref{DefinitionJS}).

The functor $\mathfrak{F}  \mapsto \mathfrak{F}(X)$ establishes an equivalence of categories:
$$\begin{CD}
\mathrm{\bf Nat}^k_n @= \mathrm{\bf Nat}^k_X
\end{CD}$$
\end{proposition}

\proof Fix any point $x_0\in X$.

The universal bundle can be simultaneously defined  for all $n$-manifolds: given an $n$-manifold $Z$, let $\mathcal{U}(Z)^k$ be the manifold formed by jets  $j^k_z\sigma$ of analytic 
isomorphisms $\, \sigma\colon X \supset U \rightarrow \, V \subset Z\,$ defined between open sets, such that $x_0 \in U$. The universal bundle for $Z$ is the projection $\,\mathcal{U}(Z)^k\to Z$, $j^k_{x_0}\sigma\mapsto \sigma (x_0)$.

As a consequence, the ``associated bundle'' functor is defined for any $n$-manifold. This allows to extend, {\it mutatis mutandis}, the arguments in the proof of Theorem \ref{TeoremaEquivalencia} to obtain an equivalence of categories:
$$\begin{CD}
  G^k_{x_0}\mathrm{\bf \text{-}Man} @= \mathrm{\bf Nat}_n^k\\
 F_0 & \longmapsto & \mathfrak{F}
\end{CD}
\qquad ,\qquad \mathfrak{F}(Z):=(\mathcal{U}^k_{x_0} (Z) \times F_0)/G^k_{x_0}\ .$$

Combining it with Theorem \ref{TeoremaEquivalencia}, it readily follows the equivalence $\,\mathrm{\bf Nat}_n^k=\mathrm{\bf Nat}_X^k$.

\qed

\subsection{Natural differential operators}

Let $F , \bar{F} \to X$ be natural bundles.

\begin{definition}\label{DiffOperatNatu}
A natural holomorphic differential operator $P \colon F_U \rightsquigarrow \bar F_U$ of finite order $\leq r$ on an open set $U \subset X$ is a morphism of natural bundles:
$$ P \colon J^r F_U \to \bar{F}_U \ . $$
\end{definition}

If $F\to X$ is a natural bundle of order $ k$, then $J^rF\to X$ is a natural bundle of order $k+r$. Therefore, Corollary  \ref{MorfismosNaturalesK} implies:

\begin{proposition}\label{Terng} 
Let $F,\bar F\longrightarrow X$ be natural bundles of order $k$.
Given any point $x_0\,$ in an open set $U$, the assignment $P \mapsto P_{x_0}$ defines a bijection:
$$\begin{CD}
\left\{\begin{aligned}
&\mathrm{Natural\ holomorphic\ operators \ } \\
& \hskip .4cm P \colon F_U \rightsquigarrow \bar F_U \ \ \mathrm{\ of\ order\ } \leq r \
\end{aligned}
\right\}
@=
\left\{\begin{aligned}
&G^{k+r}_{x_0}\text{--}\,\mathrm{equivariant\ analytic\ maps\ }\\
& \hskip 1.2cm P_{x_0} \colon J^r_{x_0} F\to\bar F_{x_0}
\end{aligned}
\right\}
\end{CD} \ . $$

In particular, any natural holomorphic differential operator is globally defined.
\end{proposition}

\medskip

\section{Natural differential operators on holomorphic forms}

Let $\,X\,$ be a compelx analytic manifold of dimension $\,n$ and let $\,T^*X$ denote the holomorphic cotangent bundle.

At any point $\,x_0 \in X$, the space of jets of holomorphic $p$-forms is isomorphic to the space of jets of $p$-forms at the origin of $\,\mathbb{C}^n$:
$$ J^r_{x_0} \left( \Lambda^p T^* X \right) \, \simeq \, J^r_0 \left( \Lambda^p T^* \mathbb{C}^n \right) \ . $$

Let us also introduce the following notations for the exterior and symmetric powers of the holomorphic cotangent space of $\,\mathbb{C}^n\,$ at the origin:
$$\,\Lambda^p:=\Lambda^p(T^*_0\mathbb{C}^n)\,\quad \mbox{ and } \quad \, S^p:=S^p(T^*_0\mathbb{C}^n) \ .$$

There exists an analytic isomorphism:
$$\begin{CD}
J^r_0(\Lambda^pT^*\mathbb{C}^n) @= (\Lambda^p)\oplus(S^1\otimes \Lambda^p)\oplus\cdots\oplus(S^r\otimes \Lambda^p)\\
j^r_0\omega & \longmapsto & (\omega^0,\omega^1,\dots  ,\omega^r)\qquad\qquad
\end{CD}$$
where $\,\omega^s\,$ is the homogeneous component of degree $\,s\,$ in the Taylor expansion (on cartesian coordinates) at the origin of the holomorphic $p$-form $\,\omega$.

If we write $\, \nabla\,$ for the standard (flat) covariant derivative of $\,\mathbb{C}^n$, then
$$\omega^s \, = \, (\nabla^s \omega)_{z=0} \ . $$

This isomorphism is not invariant under arbitrary changes of coordinates, but it is so under linear changes of coordinates. In other words, this isomorphism is equivariant with respect to the action of the linear group $\,\mathrm{GL}_n(\mathbb{C})\simeq G^1_0 \subseteq G^{r+1}_0$, although it is not equivariant with respect to the whole group $\,G^{r+1}_0$.

Proposition \ref{Terng}, combined with the previous isomorphism, imply the following:

\medskip

\begin{lemma}\label{2.1} Let $\,X\,$ be a complex analytic manifold of dimension $\,n$. There exists an injective map:
$$\begin{CD}
\left\{\begin{aligned}\mathrm{Natural\ differential\ operators\ of \ order } \leq r\\
P\colon\Lambda^{p} T^* X \,\rightsquigarrow\,\Lambda^q T^*X \qquad\qquad\end{aligned}\right\}\\
|\bigcap^{\phantom{\bigcap}}
_{\phantom{\bigcap}} \\
 \left\{\begin{aligned} &\mathrm{Analytic\ }\mathrm{GL}_n(\mathbb{C})\text{--}\,\mathrm{equivariant \ maps \ }\quad\\ \tilde P\colon  &� (S^0\otimes \Lambda^{p})\oplus\cdots\oplus(S^r\otimes \Lambda^{p})\,\longrightarrow\, \Lambda^q \, \end{aligned}\right\}
\end{CD}$$
where the relation between $\,P\,$ and $\,\tilde P\,$ is determined by the equality
\begin{align*}
P(j^r_x \omega) \, &=\,\tilde P(\omega^0,\dots,\omega^r ) \ .
\end{align*}
\end{lemma}

Later on, we will prove that this inclusion is in fact a linear isomorphism.

\subsection{Homogeneity}

Equivariance with respect to homotheties produces the following homogeneity condition:

\begin{lemma}\label{2.2} Any $\,\mathrm{GL}_n(\mathbb{C})$-equivariant map
$$\begin{CD} ( S^0\otimes \Lambda^{p})\oplus\cdots\oplus(S^r\otimes \Lambda^{p}) @>{\tilde P}>> \Lambda^q
\end{CD}$$
satisfies, for any complex number $\, \lambda\neq 0$, the following homogeneity condition
\begin{equation}\label{CondicionHomog}
\tilde P\left(\lambda^{p+0}\omega^0,\dots,\lambda^{p+r}\omega^r \right)\,=\,\lambda^q\tilde P\left(\omega^0,\dots,\omega^r \right) \ .
\end{equation}
\end{lemma}

\proof Let $\,\tau_\lambda\colon\mathbb{C}^n\to\mathbb{C}^n$, $\tau_\lambda(z)=\lambda z$, be the homothety of ratio $\,\lambda\neq 0$.

On the one hand, for any tensor $\,\omega^j\in S^j\otimes\Lambda^{p}$, covariant of order $\,p+j$, it holds $\,\tau_\lambda^*\omega^j=\lambda^{p+j}\omega^j$.

On the other, as the map $\, \tilde P\,$ is $\,\mathrm{GL}_n(\mathbb{C})$-equivariant, it holds $\,\tilde P\circ\tau_\lambda^*=\tau_\lambda^*\circ\tilde P$. Consequently:
\begin{align*}
\tilde P \left(\lambda^{p+0}\omega^0,\dots,\lambda^{p+r}\omega^r \right)\, &=\,
\tilde P\left( \tau_\lambda^*\omega^0,\dots,\tau_\lambda^*\omega^r \right) \\
& =\,\tau_\lambda^*\tilde P\left( \omega^0,\dots,\omega^r \right)  \,=\,\lambda^q\tilde P\left( \omega^0,\dots,\omega^r \right)\ .
\end{align*}

\qed

Any holomorphic map between vector spaces satisfying such a homogeneity condition has necessarily to be a polynomial, in virtue of the following elementary result  (e. g. \cite{KMS}, Thm. 24.1):

\medskip
\noindent {\bf Homogeneous Function Theorem:} {\it Let $E_1,\dots, E_k$ be finite dimensional $\mathbb{C}$-vector spaces.

Let $\, f \, \colon  \prod E_i \to \mathbb{C}$ be a holomorphic function such that there exist natural numbers $\, a_i ,  w \in \mathbb{N}\,$ satisfying:
\begin{equation}\label{CondicionHomogeneidadLemma}
 f ( \lambda^{a_1} e_1 , \ldots , \lambda^{a_k} e_k) = \lambda^w \, f(e_1 , \ldots , e_k)
\end{equation}
 for any complex number $\lambda \neq 0$ and any vector
$(e_1 , \ldots , e_k) \in \prod E_i$.

Then $f$ is a sum of monomials of
degree $(d_1,\dots,d_k)$ in the variables $e_1,\dots,e_k$ satisfying the relation
\begin{equation}\label{CondicionMonomios}
 a_1 d_1 + \cdots + a_k d_k = w \ .
\end{equation}

If there are no natural numbers  $d_1,\dots,d_k \in \mathbb{N} \cup \{ 0 \}$ satisfying this equation, then $f$ is the zero map. }
\medskip

In other words, for any finite dimensional vector space $W$, there exists a $\mathbb{C}$-linear isomorphism:
$$\begin{CD}
\left[ \text{Holomorphic maps }\, f \colon  \prod E_i \to W\,\text{ satisfying }  (\ref{CondicionHomogeneidadLemma})\right]  \\
@| \\
\bigoplus \limits _{(d_1 , \ldots , d_k)}\text{Hom}_{\mathbb{C}}(S^{d_1} E_1 \otimes \ldots
\otimes S^{d_k} E_k,\, W)
\end{CD}$$
where $(d_1, \ldots , d_k)$ runs over the non-negative integers solutions of equation (\ref{CondicionMonomios}).

A holomorphic map $\,f\colon\prod E_i\to W$, satisfying (\ref{CondicionHomogeneidadLemma}), and the corresponding linear map $\,\oplus f^{d_1\dots d_k}\in \bigoplus_{(d_1 , \ldots , d_k)}\text{Hom}_{\mathbb{C}}(S^{d_1} E_1 \otimes \ldots
\otimes S^{d_k} E_k,\, W)$, are related by the equality
$$f(e_1,\dots,e_k)\,=\,\sum_{(d_1 , \ldots , d_k)}f^{d_1\dots d_k}\left((e_1\otimes\overset{d_1}{\dots}\otimes e_1)  \otimes\cdots\cdots\otimes (e_k\otimes\overset{d_k}{\dots}\otimes e_k) \right)\ .$$

\medskip

\subsection{Invariant theory of $\mathrm{GL}_n(\mathbb{C})$}

Let $\,V\,$ be a $\mathbb{C}$-vector space of finite dimension $\,n$, and let $\,\mathrm{GL}_n(\mathbb{C})\,$ be the complex Lie group of its $\mathbb{C}$-linear automorphisms.

As any linear representation of $\,\mathrm{GL}_n(\mathbb{C})\,$ decomposes as a direct sum of irreducible representations, the following holds:

\begin{proposition}\label{Semisimple}
Let $\,E\,$ and $\,F\,$ be linear representations of $\,\mathrm{GL}_n(\mathbb{C})$, and let $\,E'\subset E\,$ be a sub-representation.
Any equivariant linear map $\,E' \to F\,$ is the restriction of an equivariant linear map $\,E \to F$.
\end{proposition}

On the other hand, the Main Theorem of the invariant theory for the general linear group states that the only $\,\mathrm{GL}_n(\mathbb{C})$-equivariant linear maps $\,\otimes^p V \longrightarrow \otimes^p V\,$ are the linear combinations of permutations of indices (\cite{LieGroups}).

In this paper, we will only use the following consequence:

\begin{proposition}\label{CorolarioMain}  The only $\,\mathrm{GL}_n(\mathbb{C})$-equivariant linear maps $\, \otimes^p V \longrightarrow \Lambda^p V\,$ are the multiples of the skew-symmetrisation operator.
\end{proposition}

Combining the previous two propositions, it follows:

\begin{corollary}\label{0.3}  Let $\,E\subseteq\,\otimes^pV\,$ be a $\,\mathrm{GL}_n(\mathbb{C})$-sub-representation. The only $\,\mathrm{GL}_n(\mathbb{C})$-equi\-va\-riant linear maps $\,E \longrightarrow \Lambda^p V\,$ are the multiples of the skew-symme\-trisation operator.
\end{corollary}

\begin{remark}\label{0.4} We will make use of the following properties of the skew-sy\-mme\-trisa\-tion operator $\,h$:\smallskip

-- $h(T\otimes T')=h(T)\wedge h(T')\,$ for any covariant tensors $\,T,T'$,\smallskip

-- $h(\omega)=q!\,\omega$ for any $q$-form $\,\omega$,\smallskip

-- $h(\nabla \omega)=q!\,\d \omega$ for any holomorphic $q$-form $\omega$ on $\,\mathbb{C}^n$, where $\,\nabla\,$ denotes the standard flat connection on the affine space.
\end{remark}

\subsection{Computation}

\begin{definition}\label{AntiConmutativos} Let us fix a finite sequence of positive integers $\,(p_1,\dots,p_m)$. Let us denote $\,\mathbb{C}\{ u_1,\dots,u_m\}\,$ the {\it anti-commutative} algebra of polynomials with complex coefficients in the variables $\,u_1,\dots, u_m$, where each variable $\,u_i\,$ is assigned degree $\,p_i$. The anti-commutative character of this algebra is expressed by the relations
$$u_iu_j\,=\, (-1)^{p_ip_j}u_ju_i\ .$$

The degree of a monomial $\,u_1^{a_1}\dots u_m^{a_m}\,$ is defined as $\,\sum a_ip_i$.

A polynomial $\,\mathbf{P}(u_1,\dots, u_m)\in \mathbb{C}\{ u_1,\dots,u_m\}\,$  is said homogeneous of degree $\,q\,$ if it is a linear combination of monomials of degree $\,q$.
\end{definition}

Let $\,\mathbf{P}(x_1,\dots, x_m)\in \mathbb{C}\{ u_1,\dots,u_m\}\,$ be a homogeneous polynomial of degree $\,q$, and let $\,\omega_1,\dots,\omega_m\,$ differential forms of degree $\,p_1,\dots,p_m\,$ on a complex manifold $\,X$. Then $\,\mathbf{P}(\omega_1,\dots,\omega_m)$, where the product of variables is replaced by the exterior product of forms, is a differential form of degree $\,q\,$ on $\,X$.

\begin{theorem}\label{Main} Let $\,p\,$ be a positive integer,  and let $\,\mathbb{C}\{ u,v\}\,$ be the anti-commutative algebra of polynomials on the variables $\,u,v$, of degree $\,\mathrm{deg}\, u=p$, $\,\mathrm{deg}\, v=p+1$.

Any natural differential operator of order $\leq r$
$$ P\colon \Lambda^{p} T^*X \  \rightsquigarrow\ \Lambda^q T^*X \ ,  $$
where $0\leq q\leq n= \dim X$, can be written as
$$P\left(j^r_x \omega \right)\,=\, \mathbf{P}(\omega,\mathrm{d}\omega )_x \quad , $$
for a unique homogeneous polynomial $\,\mathbf{P}(u,v)\in\mathbb{C}\{ u ,v \}\,$ of degree $\,q$.
\end{theorem}

\proof Due to Lemma \ref{2.1}, any such a differential operator $\,P\,$ is determined by a $\mathrm{GL}_n(\mathbb{C})$-equivariant analytic map
$$\tilde P\colon (S^0\otimes \Lambda^{p})\oplus\cdots\oplus(S^r\otimes \Lambda^{p}) \,\longrightarrow\, \Lambda^q\ ,$$
satisfying the homogeneity condition (\ref{CondicionHomog}).

The Homogeneous Function Theorem then assures that $\,\tilde P\,$ has to be a polynomial; that is to say, it is a sum $\,\oplus_{\{d_{j}\}}\, \tilde P^{\{d_{j}\}}\,$ of $\mathrm{GL}_n(\mathbb{C})$-equivariant linear maps:
$$ \tilde P^{\{ d_{j}\}}:  S^{d_{0}} (S^0\otimes\Lambda^{p}) \otimes\cdots\otimes S^{d_{r}} (S^r\otimes\Lambda^{p}) \, \longrightarrow \, \Lambda^{q} $$ where each sequence $\,\{ d_{j}\}\,$ of non-negative integers is a solution to the equation:
\begin{equation*}
p d_{0}+\cdots+(p+r)d_{r} \, = \, q \ .
\end{equation*}

Observe that this condition implies that $\,\tilde P^{\{ d_{j}\}}\,$ is in fact defined on a vector subspace of $\,\otimes^qT_0^*\mathbb{C}^n$.
Then, by Proposition \ref{0.3}, the linear map $\,\tilde P^{\{ d_{j}\}}\,$ is a multiple of the skew-symmetrisation operator $\,h$,
$$\tilde P^{\{ d_{j}\}}\,=\, \lambda^{\{ d_{j}\}}\cdot h\ .$$

The skew-symmetrisation of two symmetric indices vanishes, so that we may assume, from now on, that
$$
d_{2} =d_{3} =\ldots =d_{r}= 0 \quad  \ .$$

That is to say, we only consider solutions $\,\{ d_{0},d_{1} \}\,$ to the equation:
\begin{equation*}
p_1d_{0}+(p+1)d_{1} \, = \, q \ .
\end{equation*}

Bringing all this together,
$$P(j^r_0 \omega )\,=\,\tilde P\left((\omega)_0,(\nabla \omega )_0 \right)
$$
$$=\,\sum_{\{ d_{0}, d_1\}}\tilde P^{\{ d_{j}\}}\left(\left[(\omega)_0\otimes\overset{d_{0}}{\dots}\otimes(\omega)_0\right]\otimes\left[(\nabla \omega)_0\otimes\overset{d_{1}}{\dots}\otimes(\nabla \omega)_0\right] \right) $$
$$=\,\sum_{\{ d_{0} , d_1\}}\lambda^{\{ d_{0}, d_1\}}h\left(\left[(\omega)_0\otimes\overset{d_{0}}{\dots}\otimes(\omega)_0\right]\otimes\left[(\nabla \omega)_0\otimes\overset{d_{1}}{\dots}\otimes(\nabla \omega)_0\right]\right) $$
(using Properties \ref{0.4} of the skew-symmetrisation operator)
$$=\,\sum_{\{ d_{0} , d_1\}}\mu^{\{ d_{0} , d_1\}} \left[(\omega)_0\wedge\overset{d_{0}}{\dots}\wedge(\omega)_0\right]\wedge\left[(\d\omega)_0\wedge\overset{d_{1}}{\dots}\wedge(\d\omega)_0\right] $$
where $\mu^{\{ d_{0} , d_1\}}:=\lambda^{\{ d_{0} , d_1\}}(p!)^{d_{0}} ((p+1)!)^{d_{1}}$.

Therefore, taking
$$\mathbf{P}(u,v ):=\,\sum_{\{ d_{0}, d_1\}}\mu^{\{ d_{0}, d_1\}}u^{d_{0}}v^{d_{1}}$$
it follows
$$P(j^r_0 \omega ) \,=\, \mathbf{P}(\omega,\mathrm{d}\omega )_{x=0}\ .$$

By naturalness, we conclude
\begin{equation}\label{CondUnicidad}
P(j^r_x\omega)\,=\, \mathbf{P}(\omega,\mathrm{d}\omega )_x
\end{equation}
for any jet $\,j^r_x \omega\,$ at any point $\,x \in X$.

Finally, the uniqueness of the polynomial $\,\mathbf{P}\,$ is proved applying Lemma \ref{2.4} below to the difference of two polynomials satisfying (\ref{CondUnicidad}).

\qed

\begin{lemma}\label{2.4}
Let $\,\mathbb{C}\{ u,v\}\,$ be as in Theorem \ref{Main}.

Let $\,\mathbf{P}\in\mathbb{C}\{ u,v \}\,$ be a non-zero homogeneous polynomial of degree $\,q\leq n$. Then, the natural differential operator of order $\leq r$
$$ P\colon \Lambda^{p} T^*\mathbb{C}^n \  \rightsquigarrow \ \Lambda^q T^*\mathbb{C}^n \quad ,\quad P(j^r_x\omega )=\mathbf{P}(\omega ,\d\omega )_x  $$
is not identically zero.
\end{lemma}

\proof Up to a scalar factor, $\,\mathbf{P}\,$ is a monomial of one of the following four types, depending on the parities of $\,p\,$ and $\,q$,
$$\mathbf{P}(u,v)\, = \left\{\begin{aligned}
&\ u^0v^s\qquad (p\text{ odd, }q \text{ even}) \\
&\ u^1v^s\qquad (p\text{ odd, }q\text{ odd}) \\
&\ u^sv^0\qquad (p\text{ even, }q \text{ even}) \\
&\ u^sv^1\qquad (p\text{ even, }q\text{ odd})\end{aligned}\right.$$
where $\,s\,$ has to be taken in each case, for the monomial to have degree $\,q$.

Depending on the case, let us consider the following $\,p$-forms, which involve $q$ variables each:
$$\omega\,=\left\{\begin{aligned}
\sum_{j=1}^s z_{j , 0}\d z_{j , 1}\wedge\dots\wedge\d z_{j, {p}} \quad , \quad
\text{case }u^0v^s\\
\d z_1\wedge\dots\wedge\d z_{p }+\sum_{j=1}^s z_{j , 0}\d z_{j , 1}\wedge\dots\wedge\d z_{j , {p}}\quad , \quad \text{case }u^1v^s\\
\sum_{j=1}^s\d z_{j , 1}\wedge\dots\wedge\d z_{j ,{p}} \quad , \quad
\text{case }u^sv^0\\
z_0\d z_1\wedge\dots\wedge\d z_{p}+\sum_{j=1}^s \d z_{j , 1}\wedge\dots\wedge\d z_{j, {p}} \quad , \quad \text{case }u^sv^1
\end{aligned}\right.$$

In each of the four cases, it holds
$$P(j^r\omega)\,=\,\mathbf{P}(\omega,\d\omega)\,=\, s! \, \d z_1 \wedge \ldots \wedge \d z_q \ . $$

Hence, the differential operator $\,P\colon \Lambda^{p} T^*\mathbb{C}^n\rightsquigarrow \Lambda^q T^*\mathbb{C}^n\,$ is not identically zero.

\qed

\medskip
For the sake of completeness, let us remark that the following statement also holds, with similar reasonings in its proof:

\begin{theorem}\label{General} Let $\,p_1,\dots,p_m\,$ be positive integers,  and let $\,\mathbb{C}\{ u_1,v_1,\dots,u_m,v_m\}\,$ be the anti-commutative algebra of polynomials on the variables $\,u_1,v_1,\dots,u_m,v_m$, of degree $\,\mathrm{deg}\, u_i=p_i$, $\,\mathrm{deg}\, v_i=p_i+1$.

Any natural differential operator of order $\leq r$
$$ P\colon \Lambda^{p_1} T^*X \oplus \ldots \oplus \Lambda^{p_k} T^*X \  \rightsquigarrow\ \Lambda^q T^*X  $$
where $0\leq q\leq n = \dim X$, can be written as
$$P\left(j^r_x \omega_1,\dots,j^r_x \omega_r \right)\,=\, \mathbf{P}(\omega_1,\mathrm{d}\omega_1,\dots,\omega_k,\mathrm{d}\omega_k)_x \quad , $$
for a unique homogeneous polynomial $\,\mathbf{P}(u_1,v_1,\dots,u_k,v_k)\in\mathbb{C}\{ u_1,v_1,\dots,u_k,v_k\}\,$ of degree $\,q$.
\end{theorem}

\medskip
\noindent {\bf Acknowledgements:} We acknowledge J. A. Navarro and J. B. Sancho, who explained the theory of smooth natural operations to the first two authors, and A. Chirvasitu for raising interesting questions.

The first author has been partially supported by ICMAT Severo Ochoa project
SEV-2015-0554 (MINECO), MTM2013-42135-P (MINECO) and ANR-12-BS01-0002 (ANR).
The second author has been partially supported by Junta de Extremadura and FEDER funds. The third author has been supported by the research grant MTM2013-45935-P, Ministerio de Econom\'{\i}a y Competitividad, Spain.

\end{document}